\documentstyle{amlts}
\begin{document}
\annalsline{155}{2002}
\received{January 29, 2001}
\startingpage{281}
\def\bye{\end{document}}
 \font\tenrm=cmr10

\input boxedeps.tex 
\SetepsfEPSFSpecial 
\HideDisplacementBoxes
\def\figin#1#2{ 
$$
 {\BoxedEPSF{#1.eps scaled
#2}%
}%
$$
\noindent}
\catcode`\@=11
\font\twelvemsb=msbm10 scaled 1100
\font\tenmsb=msbm10
\font\ninemsb=msbm10 scaled 800
\newfam\msbfam
\textfont\msbfam=\twelvemsb  \scriptfont\msbfam=\ninemsb
  \scriptscriptfont\msbfam=\ninemsb
\def\msb@{\hexnumber@\msbfam}
\def\Bbb{\relax\ifmmode\let\next\Bbb@\else
 \def\next{\errmessage{Use \string\Bbb\space only in math
mode}}\fi\next}
\def\Bbb@#1{{\Bbb@@{#1}}}
\def\Bbb@@#1{\fam\msbfam#1}
\catcode`\@=12

 \catcode`\@=11
\font\twelveeuf=eufm10 scaled 1100
\font\teneuf=eufm10
\font\nineeuf=eufm7 scaled 1100
\newfam\euffam
\textfont\euffam=\twelveeuf  \scriptfont\euffam=\teneuf
  \scriptscriptfont\euffam=\nineeuf
\def\euf@{\hexnumber@\euffam}
\def\frak{\relax\ifmmode\let\next\frak@\else
 \def\next{\errmessage{Use \string\frak\space only in math
mode}}\fi\next}
\def\frak@#1{{\frak@@{#1}}}
\def\frak@@#1{\fam\euffam#1}
\catcode`\@=12



\newcommand{\R}{{\Bbb R}}
\renewcommand{\S}{{\Bbb S}}
\renewcommand{\tilde}{\widetilde}
\newcommand{\la}{\langle}	   	
\newcommand{\ra}{\rangle}          	

\renewcommand{\(}{\left(}
\renewcommand{\)}{\right)}
\renewcommand{\[}{\left[}
\renewcommand{\]}{\right]}
\newcommand{\cn}{\colon}	 
\newcommand{\ol}{\overline}   
\renewcommand{\mid}{:}
\renewcommand{\phi}{\varphi}

\def\inte{\mathop{\rm int}}
\def\conv{\mathop{\rm conv}}
\def\lan{\langle}
\def\ran{\rangle}
\def\Hess{{\rm Hess}}
\def\dist{{\rm dist}}
\def\Sym{{\rm Sym}}

\title{Shadows and convexity
                                      of  surfaces} 
\shorttitle{Shadows and convexity of surfaces} 

 \author{Mohammad Ghomi}
 \institutions{University of South Carolina 
Columbia, SC\\
{\eightpoint {\it E-mail address\/}: ghomi@math.sc.edu}\\
{\eightpoint {\it URL\/}: www.math.sc.edu/$\sim$ghomi}}

\centerline{\bf Abstract}
\vglue12pt
We study the geometry and topology of  immersed surfaces in Euclidean 3-space whose Gauss map
satisfies a certain two-piece-property, and solve the ``shadow problem" formulated by H. Wente.

\section{Introduction}

Let $M$ be a closed oriented $2$-dimensional manifold, $f\colon M\to\R^3$ be a smooth
immersion into Euclidean 3-space, and
  $n\colon M\to \S^2$ be a unit normal vectorfield, or the Gauss map, induced by $f$. Then for every
unit vector 
$u\in \S^2$  (corresponding to the direction of light)  the {\it shadow}, $S_u$, is defined by
$$
S_u:=\{\,p\in M\mid \langle n(p), u \rangle>0\,\},
$$
where $\la\cdot,\cdot\ra$ is the standard innerproduct.
If $f$ is a {\it convex embedding}, i.e., $f$ maps $M$ homeomorphically to the boundary of a
convex body,  then it is intuitively clear that
$S_u$ is a connected subset of $M$ for each $u$. In 1978, motivated by problems concerning the
stability of constant mean curvature surfaces, H. Wente
\cite{wente:email}  appears to have been the first person to study the converse
of this  phenomenon, which has since become known as the
\lq\lq shadow problem\rq\rq~\cite{mccuan:email}: {\it Does  connectedness of the shadows
$S_u$ imply that $f$ is a convex embedding}? In this paper we prove:

\proclaim{Theorem}\label{thm:1}
$f$ is a convex embedding if and only if{\rm ,} for every $u\in\S^2$,  $S_u$ is {\rm simply}
connected{\rm .}
\endproclaim 

Furthermore we show that the additional condition implied by the word ``simply" in the
above theorem is  necessary:
\vfill
\footnoterule{\ninepoint {\nineit Key words and phrases}. Shadow, skew loop, tantrix, constant mean curvature, 
two-piece-property, tight immersion.}
\vglue5pt
\eject

\proclaim{Theorem}\label{thm:2} 
There exists a smooth embedding of the torus{\rm ,} $f\colon \S^1\times\S^1\break\to\R^3${\rm ,} such that
for all $u\in \S^2$, $S_u$ is connected{\rm .}
\endproclaim 

Thus,  connectedness of the shadows  in general is
not strong enough to  ensure convexity or even determine the topology; however, we can show:

\proclaim{Theorem}\label{thm:3}
If $M$ is  topologically a sphere{\rm ,} and{\rm ,} for every $u\in\S^2${\rm ,} $S_u$ is  connected{\rm ,} then $f$ must
be a convex embedding{\rm .}
\endproclaim 

In short, the answer to the above question is yes, provided that either the shadows are
{\it simply} connected, or  $M$ is a sphere; otherwise, the answer is no. This settles
Wente's shadow problem in $3$-space. See \cite{ghomi:solution} and \cite{choe:index} for  
motivations behind this problem and  relations to  constant mean
curvature surfaces.

\numbereddemo{Note} 
The immersion $f\colon M\to\S^2$ has connected shadows if and only if for every great circle
$C\subset\S^2$, $n^{-1}(\S^2-C)$ has exactly two components. That is, the Gauss map
satisfies a 
{\it two-piece-property} \cite{cecil&ryan:book} similar to that formulated by  T. Banchoff
\cite{banchoff:tpp}, and further developed by N. Kuiper
\cite{kuiper:tight}.
\enddemo

\numbereddemo{Note} \label{note:choe}
For a great circle $C\subset\S^2$, the number of components of $n^{-1}(\S^2-C)$ has been called
the
{\it vision number} with respect to a direction  perpendicular to $C$. This terminology is due to
J. Choe, who conjectured\break \cite[p. 210]{choe:index} that there always exists a direction with respect
to which the vision number of $f\colon M\to\R^3$ is greater than or equal to $4-\chi(M)$ where $\chi$
is the Euler characteristic. Theorem \ref{thm:2} gives a counterexample to this conjecture.
\enddemo 

\section{Regularity of horizons and  shadow boundaries}

First we need to establish some basic regularity
results regarding the generic behavior of shadows.
For each $u\in\S^2$,  define the  {\it shadow function} 
$\sigma_u\cn M\to \R$ by 
$$
\sigma_u(p):=\la n(p),u\ra.
$$
 $H_u:=\sigma_u^{-1}(0)$ is called the
{\it horizon} \cite{choe:index} in the direction $u$. It is easy to see that in general
$\partial S_u\neq H_u\neq\partial S_{-u}$, where $\partial$ denotes the boundary; however,
using Sard's theorem, we can show

\proclaim{Proposition}\label{bdry-nice}
For almost all $u\in \S^2$ {\rm (}\/in the sense of Lebesgue measure\/{\rm )}
  $H_u$ is a regular  curve{\rm .}  Thus for
these $u${\rm ,}  both $\partial S_u$ and $\partial S_{-u}$ are regular curves as well{\rm .}
Further{\rm ,} if
$H_u$ is connected{\rm ,} then  $\partial S_u=H_u=\partial S_{-u}${\rm .}
\endproclaim 

We say that $\Gamma\subset M$ is a {\it regular curve} if for each $p\in\Gamma$ there
is an open neighborhood $U$ of $p$ in $M$ and a homeomorphism $\phi\colon U\to\R^2$ such that
$\phi(U\cap\Gamma)=\R$. In particular, unless stated otherwise, a regular curve needs not be
differentiable.

\demo{Proof} 
Let $T_p M$ be the tangent plane of $M$ at $p$ which  we identify with a
subspace of $\R^3$ (by identifying $T_p M$ with $f_*(T_p M)$, and parallel translating the
elements of $f_*(T_p M)$ to the origin in $\R^3$; $f_*$ denotes the {\it differential} of $f$).
Let
$UTM:=\{(p,u): p\in M, u\in T_pM,
\|u\|=1\}$ denote  the unit tangent bundle of $M$, and $\tau$ be the mapping given by
$$
UTM\ni(p,u)\stackrel{\tau}{\longmapsto}u\in\S^2.
$$
 By Sard's theorem  almost
every $u\in \S^2$ is a regular value of $\tau$; consequently, for such $u$,  $\tau^{-1}(u)$
is a regular  curve in $UTM$.

Now let $\pi$ be the mapping defined by 
$$
UTM\ni(p,u)\stackrel{\pi}{\longmapsto}p\in M,
$$
and let $u$ be
a regular value of
$\tau$.  Note that  $\pi$ is injective on $\tau^{-1}(u)$. 
  As $\tau^{-1}(u)$ is
compact, this implies that  $\pi\colon\tau^{-1}(u)\to M$
is an embedding.  Further note that
$$
\pi(\tau^{-1}(u)) =
\{p\in M: u\in T_pM\}=\{p\in M: \la n(p),u\ra=0\}=H_u.
$$
Thus $H_u$ is a regular curve.
But then, it follows that $\partial S_{u}$ and $\partial S_{-u}$ are each open in
$H_u$, which yields that $\partial S_{u}$ and $\partial S_{-u}$ are both regular
 curves as well. Finally, since these shadow boundaries are also closed in
$H_u$, it follows that whenever
$H_u$ is connected we have  $\partial S_u=H_u=\partial S_{-u}$.
\enddemo

\numbereddemo{Note} 
Suppose that there is an open set $U\subset\S^2$, such that, for all $u\in U$, both $S_u$ and
$S_{-u}$ are simply connected. Then 
$M$ is homeomorphic to~$\S^2$; because,
by the above proposition, there exists a $u_0\in U$
such that 
$H_{u_0}$ is a regular curve.
Consequently the closures $\ol S_{u_0}$ and $\ol S_{-u_0}$ are homeomorphic to disks.
 Further, since by assumption $M-{H_{u_0}}$ is made up of a pair
of simply connected components,  $H_{u_0}$ is connected.
Thus by the above proposition $\partial
S_{-u_0}=\partial S_{u_0}$.
So $M$ is homeomorphic to a pair of disks glued
together along their boundaries.
\enddemo

By {\it smooth}  we
mean differentiable of class $C^\infty$, and for convenience we always assume that the immersion
$f\colon M\to\R^3$ is  smooth, though in this paper it is enough that $f$ be $C^3$.

\numbereddemo{Note} 
The embedding $\pi\colon\tau^{-1}(u)\to M$ in the above proposition is smooth, when $u$ is a
regular value of $\tau$. In particular, $H_u$ is smooth for almost all $u\in\S^2$. To see this
let $(p, u)\in \tau^{-1}(u)$. Then $u\in T_p M$. Let $v\in T_p M$ with $\la u,v\ra=0$. Then
$c(t):=(p,\cos(t)u+\sin(t)v)$ parametrizes  the fiber $UT_pM$ of the unit tangent bundle. Note
that
$$
\tau_{*_{(p,u)}}(c'(0))=\frac{d}{dt}\tau\big(\,p,\,\cos(t)u+\sin(t)v\,\big)\Big|_{t=0}
=v\neq 0.
$$ 
On the other hand, 
$$
T_{(p,u)}(\tau^{-1}(u))=\left\{\,X\in T_{(p,u)}(UTM)\mid \tau_{*_{(p,u)}}(X)=0\,\right\}.
$$  
Thus $c'(0)\not\in T_{(p,u)}(\tau^{-1}(u))$, which
implies that
$\tau^{-1}(u)$ is never tangent to any of the fibers $UT_pM$ of the unit tangent bundle. So
$\pi|_{\tau^{-1}(u)}$ is a smooth immersion.
\enddemo

Next we need a local regularity result for the horizons and shadow boundaries. The {\it Gaussian
curvature} $K\colon M\to\R$ is defined by $K(p):=\det(n_*(p))$. 

\proclaim{Proposition}\label{prop:localnice}
If  $K(p)\neq 0$ for some $p\in M${\rm ,} then there exists a neighborhood 
$U$ of $p$ such that for all
$u\in T_p M${\rm ,} $H_u\cap U$ is a smooth regular curve and 
$\partial S_u\cap U=H_u\cap U=\partial S_{-u}\cap U${\rm .}
\endproclaim 

\demo{Proof} 
 Since
$\det(n_{*_p})=K(p)\neq 0$, then,
 by the inverse function theorem, $n$ is a  diffeomorphism between small neighborhoods $U$ of $p$
in $M$ and  $V$ of $n(p)$ in $\S^2$. Let $\S^2_u:=\{x\in\S^2:\langle x, u\rangle>0\}$. Then 
$\partial \S^2_{u}=\partial \S^2_{-u}$ is a  regular curve. Thus, since $S_u=n^{-1}(\S^2_u)$ and
$S_{-u}=n^{-1}(\S^2_{-u})$,  the proof follows.
\enddemo

\numbereddemo{Note} 
If $K(p)=0$, then $H_u$ may not be regular for {\it all} $u\in T_p M$; however, typically $H_u$
will be regular for most $u\in T_p M$; because,
for  $u\in T_p M$, the differential of $\sigma_u$ at $p$ is given by
$$
(d\sigma_u)_p(\cdot)=\la \;\cdot\,,n_{*_p}(u)\ra.
$$
So if $n_{*_p}(u)\ne 0$, e.g., $u$ is not an {\it asymptotic direction}, then
$d\sigma_u$ is  nonzero at $p$. Consequently, by the implicit function theorem,
$\sigma_u^{-1}(\sigma_u(p))=\sigma_u^{-1}(0)=H_u$ is a smooth regular curve near $p$.
\enddemo
  
\section{Critical points of height functions}

The next set of preliminary results we need involves some basic applications of Morse theory
\cite{milnor:morse}. For every $u\in \S^2$, let the {\it height function} $h_u\colon M\to\R$,
associated to the immersion $f\colon M\to\R^3$, be defined by
$$
h_u(p):=\langle f(p), u\rangle.
$$
Recall that $p$ is a {\it critical point} of $h_u$ if the differential map $(dh_u)_p\colon T_p
M\to\R$ is zero. 
Since $(dh_u)_p(\cdot)=\la \cdot, u\ra$, it follows that $p$ is a critical point of $h_u$ if and
only if $u=\pm n(p)$.
If all of its critical points are nondegenerate, $h_u$ is a {\it Morse function}.

\proclaim{Lemma}\label{prop:morse}
{\rm (i)} $h_u$ is a Morse function if and only if  $K\neq 0$ at all critical
points of $h_u${\rm .}
{\rm (ii)} $h_u$ is a Morse function for almost all $u\in\S^2${\rm .}
{\rm (iii)} The set $U\subset\S^2$ such that for all $u\in U$ $h_u$ is a Morse function is open{\rm .}
\endproclaim

Though the above is fairly well-known (e.g., see \cite[pp. 11--12]{cecil&ryan:book}), we  include a
brief proof for completeness.

\demo{Proof} 
  If $p$ is a critical
point of $h_u$, then, as a standard computation shows, the Hessian of
$h_u$ is given by
$$
\Hess h_u(\cdot, \cdot)=\pm\langle\; \cdot\,, n_{*_p}(\cdot)\rangle.
$$
Thus $h_u$ is a Morse function if and only if at each critical point $p$, 
$K(p)=\det(n_{*_p})\neq 0$. This is equivalent to requiring that both $u$ and $-u$ be regular
values of $n$, because $p$ is a critical point of $h_u$ if and only if $u=\pm n(p)$. 
Let $U\subset \S^2$ be the set of all such values. Then, by Sard's theorem, $\S^2-U$
has measure zero.  Further, since 
$M$ is compact, and the set of critical points of $n$ is closed, it follows that the set of
critical values of $n$ is closed as well, so $U$ is open.
\enddemo

The following is implicit in a paper of Chern and Lashof \cite{chern&lashofI}.

\proclaim{Lemma}\label{Chern-Lashof}
If $f$  is not a convex embedding{\rm ,} then there
exists a Morse height function $h_u$ with at least three critical points{\rm .} \hfill\qed
\endproclaim

\demo{Proof} 
Let $\#C(h_u)$ denote the number of critical points of $h_u$. Since $p$ is a critical point of
$h_u$ if and only if $n(p)=\pm u$, we have:
$$
\int_{\S^2} \#C(h_u)\, du=\int_{\S^2}\#n^{-1}(\pm u)\,du
=2\int_M|\det(n_{*})|\,dV=2\int_M|K|\,dV.
$$
The second equality above is just an application of the area formula \cite[Thm.\
3.2.3]{federer:book}, where $dV$ denotes the volume element on $M$. Suppose that $f$ is not a
convex embedding. Then, by a well-known theorem of Chern and Lashof
\cite{chern&lashofI}, 
$$
\int_M|K|\,dV> 4\pi.
$$
  Combining the
above expressions yields a lower bound for the average number of critical points:
$$
\frac{1}{4\pi}\int_{\S^2} \#C(h_u)\, du>2.
$$
So since, by Lemma \ref{prop:morse}, $h_u$ is a Morse function for almost every
$u\in\S^2$, it follows that there exists a Morse function such that $\#C(h_u)>2$.
\enddemo

\section{Triplets on the boundaries of simply connected domains}

Here we develop some elementary topological methods whose motivation  will become more
clear  in the next section.  

\numbereddemo{Definition}\label{def:nice}
By a {\it domain} we mean a connected open subset $\Omega\subset M$.
We say $\Omega$ is  {\it adjacent}  to a triplet of points $\{p_1, p_2, p_3\}\subset M$ if
$p_i\in\partial
\Omega$.
$\Omega$ is {\it regular} near $p_i$ if there are open neighborhoods $U_i$ of $p_i$ and
homeomorphisms
$\phi_i\colon U_i\to\R^2$ which map $U_i\cap\Omega$ into the upper half-plane.
 A simple closed curve $T\subset\ol\Omega$ is a
{\it triangle} of $\Omega$ (with vertices at $\{p_1, p_2, p_3\}$) if $p_i\in T$, and $T-\{p_1,
p_2, p_3\}\subset\Omega$.
\enddemo

The following lemma, though quite elementary, is more subtle than it might at first appear
(see Note \ref{note:disk}).

\proclaim{Lemma}\label{lem:triples}
Every domain $\Omega$  adjacent to $\{p_1, p_2, p_3\}$ admits a
triangle{\rm .} Further if $\Omega$ is {\rm simply} connected and regular near $p_i${\rm ,} then any pair of
such triangles may be homotoped to each other through a family of triangles of $\Omega${\rm .}
\endproclaim

\demo{Proof} 
Since $\Omega$ is open and connected, there exists a regular arc $A_{12}\subset\Omega$ whose end
points are
$p_1$ and $p_2$. Since $A_{12}$ is regular, there exists a component $(\Omega-A_{12})^+$  of
$\Omega-A_{12}$ which contains
$p_3$ in its closure.
Let
$A_{23}\subset (\Omega-A_{12})^+$ be a regular arc with end points on $p_2$ and $p_3$.
Then, similarly, there exists a component $((\Omega-A_{12})^+ -A_{23})^+$  of $(\Omega-A_{12})^+
-A_{23}$ which contains $p_1$ in its closure. Finally, let
$A_{31}\subset((\Omega-A_{12})^+ -A_{23})^+$ be a regular arc with end points at
$p_3$ and $p_1$. The union of these three arcs, and their endpoints, gives the desired triangle. 

Now suppose that $\Omega$ is simply connected and regular near $p_i$. Let $T$ and $T'$ be a pair
of triangles of
$\Omega$, and  let $A_{12}$ and $A_{12}'$ be arcs of $T$ and $T'$ respectively which connect
$p_1$ and $p_2$. Since  $\Omega$ is regular near $p_i$, we
may homotope
$A_{12}$ (while keeping its end points fixed) by a small perturbation near $p_1$ so that
$A_{12}$ and $A_{12}'$ coincide along a segment near $p_1$. Similarly, we may assume that they
coincide near $p_2$ as well. Then it remains to homotope proper subarcs of $A_{12}$ and
$A_{12}'$ which coincide at a pair of end points in $\Omega$. Since $\Omega$ is simply connected,
these  subarcs may be homotoped to each other while keeping the end points fixed. Thus
$A_{12}$ and $A_{12}'$ are homotopic through a family of arcs of $\Omega$ with end points at
$p_1$ and $p_2$. Other arcs of
$T$ may be similarly homotoped to their counterparts in $T'$, which completes the proof.
\enddemo\numbereddemo{Note} \label{note:disk}
Without the regularity assumption near $p_i$, the second claim in the above lemma is not true in
general: Suppose for instance that $\Omega\subset\R^2$ is  an open disk of radius $1$ centered at
the origin, and with segment $[0,1)$ removed.  Set $p_1=(0,0)$, $p_2=(1/2,0)$, and $p_3=(1,0)$.
Then a triangle of
$\Omega$ which lies above the $x$-axis may not be homotoped to one lying below the $x$-axis.
\enddemo

\proclaim{Proposition}\label{prop:permutation}
For a fixed orientation of $M${\rm ,} every simply connected  domain  $\Omega$ which is adjacent to
and regular near a triple of {\rm (}\/distinct\/{\rm )} points $\{p_1, p_2, p_3\}\subset M$   uniquely determines
a permutation
$\alpha_\Omega$ of
$\{p_1, p_2, p_3\}$ such that {\rm (i)}  
if $\Omega$ and $\Omega'$ 
have a triangle in common{\rm ,} then $\alpha_\Omega=\alpha_{\Omega'}${\rm ;} and {\rm (ii)} if  $\partial
\Omega=\partial\Omega'$ is a regular curve, and $\Omega$ and $\Omega'$ are distinct{\rm ,} then
$\alpha_\Omega\neq\alpha_{\Omega'}${\rm .}
\endproclaim 

\demo{Proof} 
 By  Lemma \ref{lem:triples} there exists a triangle  $T$ of $\Omega$. $T$ bounds a simply
connected subdomain $U$ of $\Omega$. Since $M$ is oriented, $U$ inherits a preferred sense of
orientation, which in turn induces an orientation, or a sense of direction, on~$T$. This
 direction induces a permutation of $\{p_1, p_2, p_3\}$ in the obvious way: If as we move
along $T$ and pass $p_1$ we reach $p_2$ before reaching $p_3$, then we set the induced
permutation $\alpha_\Omega$ to be the cycle $(p_1, p_2, p_3)$; otherwise, the induced permutation
is the cycle $(p_1, p_3, p_2)$. It is clear that these permutations depend continuously on $T$.
Thus, since by Lemma \ref{lem:triples}, all triangles of $\Omega$ are homotopic, it follows that
$\alpha_\Omega$ does not depend on the choice of $T$ and is therefore well defined; and
furthermore, if $\Omega$ and $\Omega'$ have a triangle in common then
$\alpha_\Omega=\alpha_{\Omega'}$.

Now suppose that $\partial\Omega=\partial\Omega'$ is a regular curve, and $\Omega$ and
$\Omega'$ are distinct. Then $\Omega$ and $\Omega'$ induce opposite orientations on
$\partial\Omega$ which in turn gives rise to  distinct permutations of
$\{p_1, p_2, p_3\}$ (since $\Omega$ is simply connected, $\partial\Omega$ is connected). But by  
small perturbations,
$\partial
\Omega$ may be homotoped to a triangle of $\Omega$, just as well as it may be homotoped to a
triangle of
$\Omega'$. Thus the orientations which $\Omega$
 and $\Omega'$ induce on $\partial\Omega$ are consistent with the orientations which
$\Omega$ and $\Omega'$ induce on their own triangles respectively.
 So $\alpha_\Omega\neq\alpha_{\Omega'}$.
\enddemo
\section{Proof of Theorem \ref{thm:1}}

First we show that if $f$ is a convex embedding, then $S_u$ is simply connected for all $u\in
\S^2$. To see this let $\Pi$ be a plane perpendicular to $u$ and let $\pi\colon\R^3\to \Pi$ be the
orthogonal projection. Then $D:=\pi(f(M))$ is a convex subset of $\Pi$ with interior points. In
particular, $\inte(D)$ is homeomorphic to an open disk. Since $f(M)$ is convex and by definition
$\la n(p), u\ra>0$ for all $p\in S_u$, it is not hard to verify that
$f(S_u)$ is a graph over $\inte(D)$. Thus $\pi\circ f\colon S_u\to\inte (D)$ is a
homeomorphism. 

Now we prove the other direction: Assume that for every $u\in\S^2$, $S_u$ is simply connected;
we have to show that $f$ is  a convex embedding. The proof is by contradiction:

\proclaim{Lemma}\label{lem:uv}
If $f$ is not a convex embedding{\rm ,} then there exists a pair of orthogonal vectors $u_0${\rm ,}
$v_0\in\S^2$ such that {\rm (i)} $h_{u_0}$ is a Morse function with at least three critical points{\rm ,} and
 {\rm (ii)} 
$\partial S_{v_0}=H_{v_0}=\partial S_{-v_0}$ is a  regular  curve{\rm .}
\endproclaim

\demo{Proof} 
By Lemma \ref{Chern-Lashof}, there exists a unit vector
$u\in\S^2$ such that the corresponding height function $h_u$ is a Morse function
and has at least three critical points. Further, it follows from Lemma \ref{prop:morse},
 that this $u$ may be chosen from an open set $U\subset\S^2$.

Let $u^\perp:=\{v\in\S^2:\la u, v\ra=0\}$. Then $U^\perp:=\cup_{u\in U}u^\perp$ is
open. Consequently, by Proposition \ref{bdry-nice}, there exits a $v_0\in u_0^\perp\subset
U^\perp$ such that
$H_{v_0}$ is a regular curve. Further, since the complement of $H_{v_0}$ consists of a pair of
simply connected domains, $H_{v_0}$ is connected. Thus, again by Proposition \ref{bdry-nice},
$\partial S_{v_0}=H_{v_0}=\partial S_{-v_0}$ is a regular curve. 
\enddemo

Let $\widehat v_0\in\S^2$ be a vector orthogonal to both $u_0$ and $v_0$, and set
\begin{equation}\label{eq:v}
v(\theta):=\cos(\theta)\, v_0+ \sin(\theta)\, \widehat v_0.
\end{equation}
Let $p_i$, $i=1$, $2$, $3$, be a fixed triple of (distinct) critical points of $h_{u_0}$.

\proclaim{Lemma}\label{lem:theta}
For all
$\theta\in\R$, $S_{v(\theta)}$ is a domain adjacent to and regular near $p_i${\rm .}
\endproclaim

\demo{Proof} 
 If $p_i$ is a critical point of $h_{u_0}$, then $n(p_i)=\pm u_0$. So $\sigma_{v(\theta)}(p_i)=\la
v(\theta),\pm u_0\ra =0$, which yields that $p_i\in H_{v(\theta)}$. Since $h_{u_0}$ is a Morse
function, then, by Lemma
\ref{prop:morse},
 $K(p_i)\neq 0$. So by Proposition \ref{prop:localnice}, there exists  a neighborhood $U_i$ of
$p_i$ such that $\partial S_{v(\theta)}\cap U_i=H_{v(\theta)}\cap U_i= \partial S_{-v(\theta)}\cap
U_i$, which completes the proof.
\enddemo

It now follows from Proposition \ref{prop:permutation} that each $S_{v(\theta)}$ induces a
permutation of $\{p_1, p_2, p_3\}$ which we denote by $\alpha_\theta:=\alpha_{(S_{v(\theta)})}$.
Further, by the same proposition and since $\partial S_{v_0}=\partial S_{-v_0}$ is a regular
curve, it follows that $\alpha_0\neq\alpha_\pi$. On the other hand, 
letting $\Sym$ denote the symmetric group, we claim that
the mapping 
$$
\R\ni\theta\longmapsto\alpha_\theta\in \Sym\big(\{p_1, p_2, p_3\}\big)
$$ 
is locally constant, which, since $[0,\pi]$ is
connected, would imply that
$\alpha_0=\alpha_\pi$. This contradiction, which would complete the proof, follows from
Proposition \ref{prop:permutation} and the following:  

\proclaim{Lemma}\label{lem:continuity}
For each $\theta_0\in\R$ there exists an $\varepsilon>0$ such that if
$|\theta-\theta_0|<\varepsilon$ then $S_{v(\theta)}$ and $S_{v(\theta_0)}$ have a common triangle
{\rm (}\/with vertices at $\{p_1, p_2, p_3\}${\rm ).}
\endproclaim

\demo{Proof} 
 Recall that, since $h_{u_0}$ is a Morse
function, then, by Lemma \ref{prop:morse}, $K(p_i)\neq 0$ which yields that $n$ is a local
diffeomorphism at $p_i$. Therefore, by Proposition \ref{prop:localnice}, in a neighborhood $W$ of
$\{p_1, p_2, p_3\}$,
$\partial S_{v(\theta)}=H_{v(\theta)}=n^{-1}(v^\perp(\theta))$ where $v^\perp(\theta)$ denotes
the great circle in
$\S^2$ orthogonal to $v(\theta)$. So, since $v^\perp(\theta)$ depends continuously on $\theta$,
it follows that,  in $W$,
$\partial S_{v(\theta)}$ depends continuously on $\theta$ as well.

Let $T$ be a triangle of $S_{v(\theta_0)}$. Since $S_{v(\theta_0)}$ is open, after a perturbation
of $T$ we may assume that the arcs of $T$ are smooth and meet $\partial S_{v(\theta_0)}$
transversely (recall that, by Proposition \ref{prop:localnice}, $\partial S_{v(\theta_0)}$ is
smooth near $p_i$). Thus, by the above paragraph, it follows that if
$|\theta-\theta_0|<\varepsilon_1$, for some sufficiently small
$\varepsilon_1>0$, then
$T$ meets
$\partial S_{v(\theta)}$ transversely as well.  Then it follows that for some  neighborhood
$W$  of 
$\{p_1, p_2,p_3\}$, $(T-\{p_1, p_2,p_3\})\cap W\subset S_{v(\theta)}$
for all $\theta$ such that $|\theta-\theta_0|<\varepsilon_1$.

Next note that $T-W$ is compact, and the mapping $\theta\mapsto\sigma_{v(\theta)}$ is continuous;
therefore, since by assumption
$\sigma_{v(\theta_0)}>0$ on
$T-W$, it follows that there exists an $\varepsilon_2>0$ such that  $\sigma_{v(\theta)}>0$ on $T-W$
for all $\theta$ such that $|\theta-\theta_0|<\varepsilon_2$.
This yields that $T-W\subset S_{v(\theta)}$ for all $\theta$ such that
$|\theta-\theta_0|<\varepsilon_2$. 

From the previous two paragraphs it follows that,  setting $\varepsilon:= \min\{\varepsilon_1,
\varepsilon_2\}$, we have  $(T-\{p_1, p_2,p_3\})\subset S_{v(\theta)}$ for all $\theta$ such that
$|\theta-\theta_0|<\varepsilon$, which completes the proof.
\enddemo

\numbereddemo{Note} 
Theorem \ref{thm:1} does not remain valid if the shadows are defined as the sets where 
$\la n(p), u\ra\geq 0$. For instance, the standard torus of revolution would be a counterexample.
\enddemo

\numbereddemo{Note} 
Theorem \ref{thm:1} does not remain valid without the compactness assumption; the
hyperbolic paraboloid given by the graph of $z=xy$ would be a counterexample. This follows because
here the unit normal vectorfield
$n$ is a homeomorphism into a hemisphere. Thus the preimage of any open
hemisphere under $n$ is simply connected.
\enddemo
\vglue-6pt
\section{Proof of Theorem \ref{thm:2}}
\vglue-9pt
\numbereddemo{Definition}\label{def:skewloop}
We say an immersion $\gamma\colon\S^1\simeq\R/2\pi\to\R^3$
is a {\it skew loop} if it has no pair of distinct parallel
tangent lines; i.e, 
$$
\gamma'(t)\times\gamma'(s)\neq 0
$$
for all
$t,\,s\in[0, 2\pi)$, $t\neq s$. 
\enddemo

A specific example of a skew loop, formulated by Ralph Howard [11], is as follows:

\numbereddemo{Example}
 Let
$\gamma(t):=(x(t), y(t), z(t))$, where
\begin{eqnarray*}
x(t) &:= &-\cos(t)-\frac{1}{20}\cos(4t)+\frac{1}{10}\cos(2t), \\
y(t) &:= &+\sin(t)+\frac{1}{10}\sin(2t)+\frac{1}{20}\sin(4t),\\ 
z(t) &:= &-\frac{46}{75}\sin(3 t)-\frac{2}{15}\cos(3t)\sin(3 t) ,
\end{eqnarray*}
and $t\in[0,2\pi]$. A computation of the tangential indicatrix $T(t):=\gamma'(t)/\|\gamma'(t)\|$
shows that $T(t)\neq\pm T(s)$ for all $t$, $s\in[0,2\pi)$, $t\neq s$. Thus
$\gamma$ is a skew loop. Figure~1 shows the pictures of a
 tube built around $\gamma(\S^1)$.
 \figin{ghomitubesc}{1000}

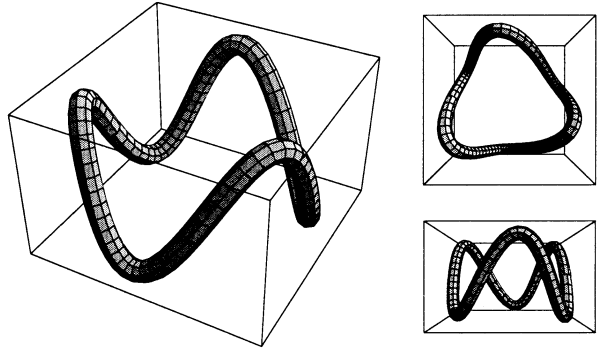
\centerline{Figure 1}
\enddemo

If $\gamma\colon\S^1\to\R^3$ is an immersion, then the unit normal bundle of $\gamma$ consists of
all pairs $(p,\nu)\in \S^1\times\S^2$ such that $\la \gamma'(p), \nu\ra=0$. 
Since this bundle  is
 homeomorphic to a torus, the following
proposition yields Theorem \ref{thm:2}.

\proclaim{Proposition}
Let $\gamma\colon\S^1\to\R^3$ be a skew loop and $M$ be the unit normal bundle of $\gamma${\rm .} For
$\varepsilon>0${\rm ,} define
$f_\varepsilon\colon M\to\R^3$ by 
$$
f_\varepsilon(p,\nu):=\gamma(p)+\varepsilon\, \nu.
$$
 Then{\rm ,}  for $\varepsilon$ sufficiently small{\rm ,}
$f_\varepsilon$ is a smooth immersion{\rm ,} and for all
$u\in\S^2${\rm ,} $S_u$ is connected{\rm .} If $\gamma$ is an embedding{\rm ,} then $f_\varepsilon$ is an embedding as
well{\rm .}
\endproclaim 

\demo{Proof} 
That $f_\varepsilon$ is a smooth immersion and is an embedding when $\gamma$ is embedded follows from
the tubular neighborhood theorem.  Let $n\colon
M\to\S^2$ be the unit normal vector field given by $n(p,\nu)=\nu$, and
$\pi\colon M\to\S^1$ be  given by $\pi(p, \nu)=p$. For every
$p\in\S^1$, let $F_p:=\pi^{-1}(p)$ be the corresponding fiber. Note that $n$
embeds  $F_p$ into the great circle in $\S^2$ which lies in the plane perpendicular to
$T(p)$. 
Further recall that
$S_u=n^{-1}(\S^2_u)$ where
$\S^2_u$ is the open hemisphere determined by $u$.
  Thus there are only two
possibilities for each $p\in\S^1$: either
$F_p$ intersects
$S_u$ in an open half-circle, or $F_p$ is  disjoint from $S_u$. The latter occurs if and
only if $T(p)$ is parallel to $u$, which, since
$\gamma$ is skew, can occur  at most once. Hence, it follows that  $S_u$ is
either homeomorphic to a disk or an annulus. In particular,
$S_u$ is connected for every $u\in\S^2$. 
\enddemo

\numbereddemo{Question}
Let $M$ be a closed oriented $2$-dimensional manifold with topological genus $g(M)\geq 2$. Does
there exist an embedding, or an immersion, $f\colon M\to\R^3$ such that $S_u$ is connected for all
$u\in\S^2$?
\enddemo

\numbereddemo{Note} 
Skew loops were first discovered by B. Segre
\cite{segre:global} to disprove a conjecture of H. Steinhaus (see also
\cite{porter:note}). More recently,
it has been shown that there exists a skew loop in each knot class
\cite{wu:knots}, and every pair of knots may be realized with the same tangential
indicatrix \cite{adams:crossing}.
\enddemo

\numbereddemo{Note} \label{note:skew}
A general procedure for constructing skew loops is as follows.
Let $T\subset\S^2$ be a smooth simple closed curve such that (i) the origin is
contained in the interior of the convex hull of $T$, $(0,0,0)\in\inte\conv T$, and
(ii) $T$ does not contain any pair of antipodal points, $T\cap -T=\emptyset$.  Figure
2 shows an example.
 \figin{tripodsc}{900}
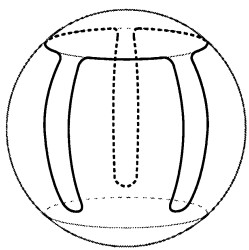
\centerline{Figure 2}
\vglue12pt
\noindent Let $T(s)$, $s\in\R$, denote a periodic parametrization of $T$ by arclength. So,
assuming
$T$ has total length $L$, we have $T(s+L)=T(s)$. Since
$(0,0,0)\in\inte\conv T$, there exists a  function $\rho(s)$ with period $L$
such that $\int_0^L\rho(s)T(s)\,ds=0$ \cite[p. 168]{gromov:PDR}. Set
$$\gamma(t):=\int_0^t \rho(s)T(s)\,ds.$$ Then $\gamma(t+L)=\gamma(t)$. Further,
$\gamma'(t)/\|\gamma'(t)\|=T(t)$. Thus $\gamma$ is a  closed curve whose
tangential spherical image coincides with $T$.  Hence $\gamma$ is a skew loop.
\enddemo

\numbereddemo{Note} 
With the sole exception
of ellipsoids, every closed surface immersed in $\R^3$ admits a skew loop \cite{ghomi&solomon}.
\enddemo

\section{Proof of Theorem \ref{thm:3}}

We follow a modified outline of the proof of Theorem \ref{thm:1}, which again proceeds by
contradiction. Suppose that $M$ is homeomorphic to $\S^2$ and $S_u$ is connected for all
$u\in \S^2$. If
$f$ is not a convex embedding, let
$u_0$ and
$v_0$ be as in Lemma \ref{lem:uv}, and  $v(\theta)$ be as defined by (\ref{eq:v}). 

\numbereddemo{Definition}
 The {\it augmented shadow} $\tilde{S}_{v(\theta)}$ is the union of
$S_{v(\theta)}$ with all components $X$ of $H_{v(\theta)}$ such that $U-X\subset
S_{v(\theta)}$ for an open neighborhood $U$ of $X$. 
\enddemo

Then $\tilde{S}_{v(\theta)}$ satisfies the conditions of the following lemma:

\proclaim{Lemma}
If $U\subset\S^2$ is a connected open set{\rm ,} and  $\S^2-U$ is also connected and has an
interior point{\rm ,} then $U$ is simply connected{\rm .}
\endproclaim

\demo{Proof} 
Let $p$ be an interior point of $\S^2-U$. Then the stereographic  projection  maps $U$ into a
connected open set with connected complement. Thus, by \cite[Thm.
11.4.1]{greene&krantz:book}, $U$ is simply connected.
\enddemo

So $\tilde{S}_{v(\theta)}$ is simply connected. Further:

\proclaim{Lemma}
For all
$\theta\in\R${\rm ,} $\tilde{S}_{v(\theta)}$ is a domain adjacent to and regular near 
$p_i${\rm .}
\endproclaim

\demo{Proof} 
This follows just as in the proof of Lemma \ref{lem:theta}, once we observe that
whenever $\partial S_{v(\theta)}=H_{v(\theta)}=\partial S_{-v(\theta)}$ is regular in some open
neighborhood, then
$\partial\tilde{S}_{v(\theta)}$, and
$\partial S_{v(\theta)}$ coincide within that neighborhood.
\enddemo

Thus each $\theta$ induces a permutation
$\tilde\alpha_\theta:=\alpha_{(\stackrel{\sim}{S}_{v(\theta)})}$ of
$\{p_1, p_2, p_3\}$ which satisfies the enumerated properties in Proposition
\ref{prop:permutation}.  In particular $\tilde\alpha_0\neq\tilde\alpha_\pi$, because since
$\partial S_{v(0)}=\partial S_{-v(0)}$ is by Lemma \ref{lem:uv} a regular curve, it follows that
$\partial
\tilde{S}_{v(0)}=\partial
\tilde{S}_{-v(0)}$ is a regular curve as well. So it remains to verify the following lemma which 
shows that $\theta\mapsto\tilde\alpha_{\theta}$ is locally constant. This would yield that 
$\tilde\alpha_0=\tilde\alpha_\pi$ which is the desired contradiction.

\proclaim{Lemma}
For each $\theta_0\in\R$ there exists an $\varepsilon>0$ such that if
$|\theta-\theta_0|<\varepsilon$ then $\tilde{S}_{v(\theta)}$ and $\tilde{S}_{v(\theta_0)}$ have a
common triangle {\rm (}\/with vertices at $\{p_1, p_2, p_3\}${\rm ).}
\endproclaim

\demo{Proof} 
This is an immediate consequence of Lemma \ref{lem:continuity} where it was proved that
$S_{v(\theta)}$ and $S_{v(\theta_0)}$ have a triangle in common (the
proof of Lemma \ref{lem:continuity} makes no use of the simply connectedness assumption on
$S_{v(\theta)}$).
\enddemo
\demo{Acknowledgements}
The author is  grateful to Ralph Howard, whose notes \cite{howard:shadow}
provided a helpful basis for the exposition of this paper. Also I 
thank John McCuan for bringing the shadow problem to my attention, and Henry Wente
for his encouragement. Finally, I would like to acknowledge the
hospitality of the Mathematics Department at the University of California at
Santa Cruz, where parts of this work were completed.


\begin{references}

\bibitem{adams:crossing}
\name{C.~Adams, C.~Lefever, J.~Othmer, S.~Pahk, A.~Stier}, and \name{J.~Tripp}, 
  An introduction to the supercrossing index of knots and the crossing
  map, preprint.

\bibitem{banchoff:tpp}
\name{T.~F. Banchoff}, 
  The two-piece property and tight $n$-manifolds-with-boundary in
  ${E}\sp{n}$, 
  {\it  Trans.\ Amer.\ Math.\ Soc.}\ {\bf 161} (1971), 259--267.

\bibitem{cecil&ryan:book}
\name{T.~E. Cecil} and \name{P.~J. Ryan}, 
  {\it  Tight and Taut Immersions of Manifolds}, 
  Pitman (Advanced Publishing Program), Boston, Mass., 1985.

\bibitem{chern&lashofI}
\name{S.-s. Chern} and \name{R.~K. Lashof}, 
  On the total curvature of immersed manifolds, 
  {\it  Amer.\ J. Math.}\ {\bf 79} (1957), 306--318.

\bibitem{choe:index}
\name{J.~Choe}, 
  Index, vision number and stability of complete minimal
surfaces,
  {\it  Arch.\ Rational Mech.\ Anal.}\ {\bf 109} (1990), 195--212.

\bibitem{federer:book}
\name{H.~Federer}, 
  {\it  Geometric Measure Theory}, 
   {\it Die Grundlehren der mathematischen Wissenschaften\/} {\bf 153},  Springer-Verlag,  New York, 1969.


\bibitem{ghomi:solution}
\name{M.~Ghomi}, 
  Solution to the shadow problem in $3$-space, 
  in {\it  Minimal Surfaces, Geometric Analysis and Symplectic
Geometry\/},
  2000, {\it Adv.\  Stud.\   Pure Math\/}.\  (To appear); preprint
  available at {www.math.sc.edu/$\sim$ghomi}.

\bibitem{ghomi&solomon}
\name{M.~Ghomi} and \name{B.~Solomon}, 
  Skew loops and quadric surfaces, 
 preprint available at\break
 {www.math.sc.edu/$\sim$ghomi}.

\bibitem{greene&krantz:book}
\name{R.~E. Greene} and \name{S.~G. Krantz}, 
  {\it  Function Theory of One Complex Variable},
  John Wiley \& Sons Inc., New York, 1997.

\bibitem{gromov:PDR}
\name{M.~Gromov}, 
  {\it  Partial Differential Relations},
  Springer-Verlag, New York,  1986.

\bibitem{howard:shadow}
\name{R.~Howard}, 
  Mohammad {G}homi's solution to the shadow problem, 
   Lecture notes, available at www.sc.edu/$\sim$howard.

\bibitem{kuiper:tight}
\name{N.~H. Kuiper}, 
  Geometry in curvature theory,
  in {\it  Tight and Taut Submanifolds\/} (Berkeley, CA, 1994),
pp.\ 
  1--50,  Cambridge Univ. Press, Cambridge, 1997.

\bibitem{mccuan:email}
\name{J.~McCuan}, 
  Personal e-mail,
  June 23, 1998.

\bibitem{milnor:morse}
\name{J.~Milnor},
  {\it  Morse Theory},
  Based on lecture notes by M.\ Spivak and R.\ Wells, {\it Ann.\
of Math.\ Studies\/} {\bf 51}, Princeton Univ.\ Press, Princeton, NJ, 1963.

\bibitem{porter:note}
\name{J.~R. Porter},
  A note on regular closed curves in ${E}\sp{3}$,
  {\it  Bull.\ Acad.\ Polon.\ Sci.\ S{\rm \'{\it e}}r.\ Sci.\ Math.\ Astronom.\ Phys.}\
  {\bf 18} (1970), 209--212.

\bibitem{segre:global}
{\it B.~Segre},
  Global differential properties of closed twisted curves,
  {\it  Rend.\ Sem.\ Mat.\ Fis.\ Milano} {\bf 38} (1968), 256--263.


\bibitem{wente:email}
\name{H.~C. Wente},
  Personal e-mail,
  January 9, 1999.

\bibitem{wu:knots}
\name{Y.-Q. Wu},
  Knots and links without parallel tangents.
  {\it  Bull.\ London Math.\ Soc\/}., to appear.

\end{references}
\end{document}